# Locating the roots of a quadratic equation in one variable through a Line-Circumference (LC) geometric construction in the plane of complex numbers

Daniel Alba Cuellar
albacd@cimat.mx

**Abstract**

This paper describes a geometrical method for finding the roots $r_1, r_2$ of a quadratic equation in one complex variable of the form $x^2 + c_1 x + c_2 = 0$, by means of a Line $L$ and a Circumference $C$ in the complex plane, constructed from known coefficients $c_1, c_2$. This Line-Circumference (LC) geometric structure contains the sought roots $r_1, r_2$ at the intersections of its component elements $L$ and $C$. Line $L$ is mapped onto Circumference $C$ by a Möbius transformation. The location and inclination angle of $L$ can be computed directly from coefficients $c_1, c_2$, while $C$ is constructed by dividing the constant term $c_2$ by each point from $L$. This paper describes the technical details for the quadratic LC method, and then shows how the quadratic LC method works through a numerical example. The quadratic LC method described here, although more elaborate than the traditional quadratic formula, can be extended to find initial approximations to the roots of polynomials in one variable of degree $n \geq 3$. As an additional feature, this paper also studies an interesting property of the rectilinear segments connecting key points in a quadratic LC structure.

**Description of the LC method for finding the roots of a quadratic equation**

Let's consider the quadratic equation

$$x^2 + c_1 x + c_2 = (x - r_1)(x - r_2) = 0. \tag{1}$$

In equation (1), we will assume that both the coefficients $c_1, c_2$ of the polynomial on the left and the roots $r_1, r_2$ in the central factorization are elements of the set of complex numbers, denoted as $\mathbb{C}$. Moreover, we will suppose that roots $r_1, r_2$ are both different from zero, and different from each other.

The geometrical method described below for finding the unknown roots $r_1, r_2$ from known coefficients $c_1, c_2$, can be used when the roots $r_1, r_2$ of equation (1) are on a line $L_1 \subset \mathbb{C}$ that does not pass through the origin (i.e., $0 \notin L_1$, $L_1$ being a continuous line of infinite length, with no gaps). This line $L_1$ can be expressed as a parametric trajectory of the form

$$L_1: p_1 + t v_{\theta^*}, \tag{2}$$

where $p_1, v_{\theta^*} \in \mathbb{C}$ are fixed values, and $t \in \mathbb{R}$ is a parameter that determines each point contained in $L_1$. Point $p_1$, defined as $p_1 := -c_1/2 = (r_1 + r_2)/2$, is called the *fixed point of $L_1$*, while point $v_{\theta^*}$, defined as $v_{\theta^*} := e^{i\theta^*} = \cos\theta^* + i\sin\theta^*$ ($i = \sqrt{-1}$), is called the *direction vector of $L_1$*. Notice that $\theta^*$ is the argument of $v_{\theta^*}$, and also, $|v_{\theta^*}| = 1$; that is to say, $v_{\theta^*}$ is a *unit direction vector of $L_1$*, and $\theta^*$ is the inclination angle of line $L_1$ that causes $L_1$ to contain $r_1, r_2$. Notice that $p_1$ is the midpoint, on $L_1$, between roots $r_1, r_2$, so $|r_1 - p_1| = |r_2 - p_1| > 0$.

How do we compute the angle $\theta^*$ from coefficients $c_1$, $c_2$ in equation (1)? From line $L_1$, it is possible to construct a semi-line of the form

$$L_d: (p_1 + tv_{\theta^*})(p_1 - tv_{\theta^*}) = p_1^2 + t^2(-v_{\theta^*}^2), \qquad (3)$$

with fixed point $p_1^2 = c_1^2/4$ and direction vector $-v_{\theta^*}^2$. Notice that the inclination angle of $L_d$ in (3) is $\arg(-v_{\theta^*}^2) = 2\theta^* + \pi$ radians, and that $r_1 r_2 = c_2 \in L_d$. This means that we can construct $L_d$ directly from the polynomial coefficients $c_1$, $c_2$, and therefore we can derive angle $\theta^*$ from the inclination angle of semi-line $L_d$.

Operationally, we obtain $\theta^*$ in the following way:

1. We construct a direction vector for $L_d$ from its two known points $p_1^2$, $c_2$; let's call this direction vector $v_d$. In this way, $v_d = c_2 - p_1^2$.
2. Since $\arg(v_d) = 2\theta^* + \pi$, we have that $\theta^* = \arg(-v_d)/2$.

In conclusion,

$$\theta^* = \arg(c_1^2/4 - c_2)/2. \qquad (4)$$

Now, we know the fixed point $p_1$ and the inclination angle $\theta^*$ of the line $L_1$ that contains the roots $r_1, r_2$ of equation (1), so we can trace it on the plane $\mathbb{C}$. To determine the location of the roots $r_1, r_2$ within $L_1$, let's consider the circumference

$$C: c_2/L_1 = c_2/(p_1 + tv_{\theta^*}). \qquad (5)$$

Notice that $C \to 0$ when $t \to \pm\infty$. Notice also that roots $r_1, r_2$ are contained in circumference $C$. Expression (5) actually comes from a bijective mapping $L_1 \to c_2/L_1$, called Möbius transformation. A demonstration that expression (5) is indeed a circumference in the plane $\mathbb{C}$ can be found in Appendix 1.

From the above, we see that the unknown roots $r_1, r_2$ are the intersections between $L_1$ and $C$.

To determine the center and radius of circumference $C$ defined by expression (5), it is sufficient to know two distinct points on it; one of them could be $w_1 = c_2/p_1$, and another could be $w_2 = c_2/(p_1 + v_{\theta^*})$. It is then possible to prove that the center $c$ of circumference $C$ is the point

$$c = \frac{-y_1(x_2^2 + y_2^2) + y_2(x_1^2 + y_1^2)}{2x_1 y_2 - 2x_2 y_1} + i \frac{x_1(x_2^2 + y_2^2) - x_2(x_1^2 + y_1^2)}{2x_1 y_2 - 2x_2 y_1}, \qquad (6)$$

where $x_1 = \text{Re}(w_1)$, $y_1 = \text{Im}(w_1)$, $x_2 = \text{Re}(w_2)$, and $y_2 = \text{Im}(w_2)$. The radius of circumference $C$ is simply $|c|$ (see Appendix 2 for details).

In this way, we have constructed a computable geometrical method to determine the location of the roots of quadratic equations of the form (1); as previously stated, this method will work as long as the distinct roots $r_1, r_2$ are contained in a continuous line $L_1 \subset \mathbb{C}\setminus\{0\}$ with no gaps.

In summary, the steps to computationally implement this geometrical method are as follows:

LC method for locating the roots of a quadratic equation

**LC method for finding the roots $r_1, r_2 \in \mathbb{C}$ of quadratic equation $x^2 + c_1 x + c_2 = 0$**

a) Compute the inclination angle $\theta^*$ for line $L_1$ by using expression (4).
b) Determine line $L_1$ defined by expression (2), using fixed point $p_1 = -c_1/2$ and direction vector $v_{\theta^*} = \cos\theta^* + i\sin\theta^* = e^{i\theta^*}$.
c) Determine circumference $C: c_2/L_1$ with center $c$ given by (6), and with radius $|c|$.
d) Locate the intersections between $L_1$ and $C$. These are the roots $r_1, r_2$ of equation (1).

**Example**

Find the roots of equation

$$x^2 + (-1 - 7i)x + (-18 + i) = 0 \tag{7}$$

using the LC method described above.

**Solution**

According to the form of equation (1), we designate coefficients in equation (7) as $c_1 = -1 - 7i$, and $c_2 = -18 + i$. Now let's carry out the steps of the LC method described above:

a) First, we compute the inclination angle $\theta^*$ for line $L_1$ via the coefficients $c_1, c_2$:
$\theta^* = \arg(c_1^2/4 - c_2)/2 = \arg(6 + 2.5i)/2 = 0.1973956$ radians $= 11°\,18'\,35.757"$.

b) Now, we obtain the parametric expression for line $L_1$:
$L_1: (1 + 7i)/2 + te^{0.1973956i} = (0.5 + 3.5i) + t(0.9805807 + 0.1961161i)$.

c) From b), we immediately obtain the parametric expression for $C$:
$C: c_2/L_1 = (-18 + i)/[(0.5 + 3.5i) + t(0.9805807 + 0.1961161i)]$.
According to (6), this parametric expression $C$ is a circumference with center $c = 0.676471 + 2.617647i$ and radius $|c| = 2.703644$. Figure 1 shows the parametric trajectories $L_1$ and $C$ obtained in this example.

d) From Figure 1 we can see that the intersections between line $L_1$ and circumference $C$ occur at points $r_1 = -2 + 3i$ and $r_2 = 3 + 4i$. We can see that these values obtained are in fact the ones that satisfy equation (7). It is perfectly feasible to construct a computational algorithm that numerically approximates the intersections between a line and a circumference in the plane of complex numbers $\mathbb{C}$, given a fixed point and a direction vector for the line, as well as the center and radius of the circumference.

LC method for locating the roots of a quadratic equation

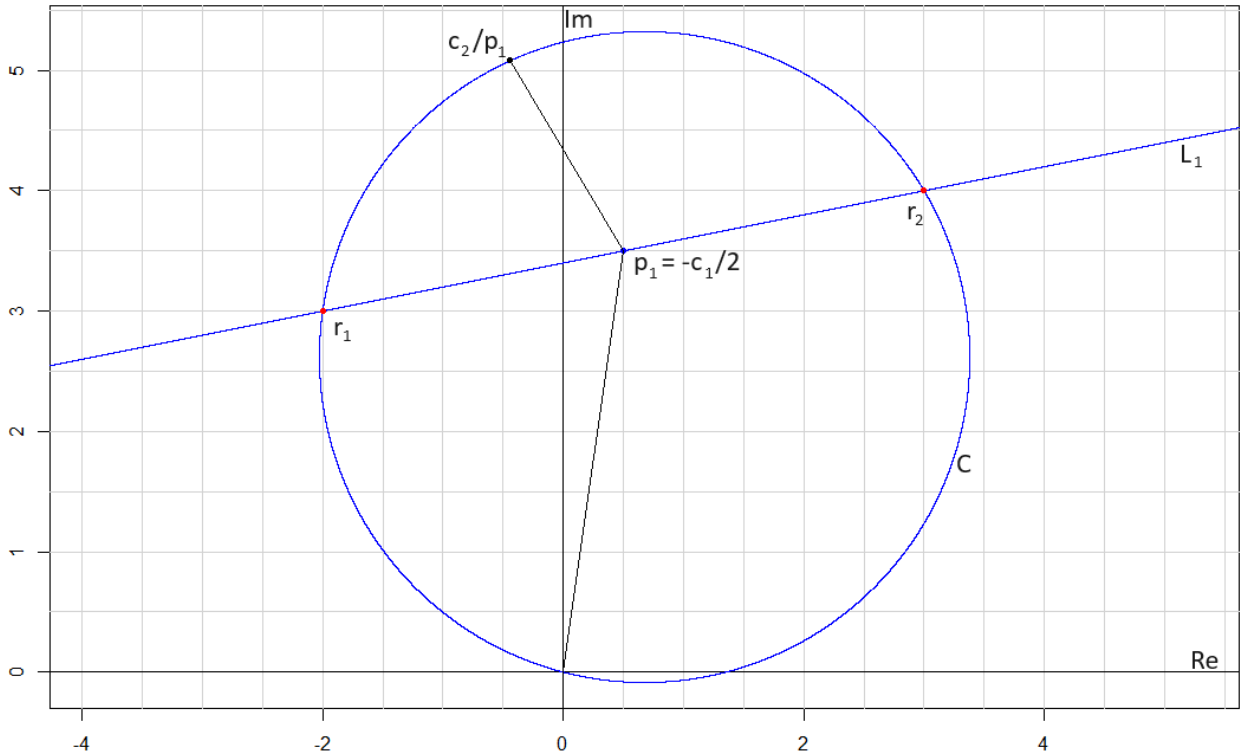

**Figure 1**. Line $L_1$ and circumference $C$ associated with equation (7), plotted on the plane $\mathbb{C}$. The intersections between these two trajectories coincide with the roots of (7).

**Conclusion**

In this paper we have described a geometrical method that helps us locate the roots of a quadratic equation in a single complex variable of the form (1), under the assumption that the roots are different from each other and are on a line in the complex plane that does not pass through the origin. Clearly, the quadratic LC method described here is more elaborate than the general formula for solving quadratic equations; however, the quadratic LC method serves as the basis for a general algorithm that can be used to find initial approximations to the roots of polynomial equations in one variable of degree $n \geq 3$. For more details, see [1].

As an additional note, referring to Figure 1, it is possible to verify numerically that angle $\angle 0 p_1 r_1$ is equal to angle $\angle r_1 p_1 (c_2/p_1)$, or equivalently, that angle $\angle 0 p_1 r_2$ is equal to the angle $\angle r_2 p_1 (c_2/p_1)$; this does not happen by chance. In Appendix 3 it is shown that this is in fact a property that holds for geometric constructions with characteristics similar to the one in Figure 1. From this, we can conclude that an alternative method for determining the inclination angle $\theta^*$ for line $L_1$ can be based on the bisection of angle $\angle 0 p_1 (c_2/p_1)$, which can be constructed directly from the coefficients $c_1$, $c_2$ of equation (1).

**Appendix 1**

**Proposition**. If $z \in \mathbb{C}$ is on a line $L$ that does not pass through the origin, then $b/z$ ($b \in \mathbb{C}, b \neq 0$) is on a circumference $C$ that passes through the origin without containing it.

**Demonstration**. In notes by Professor Carl Eberhart (University of Kentucky), available at [2], it is shown that if $z$ is on a line $L$ that does not pass through the origin, then $1/z$ is on a circumference $C$ that passes through the origin without containing it. We shall see the demonstration of this lemma below, based on the notes [2].

We start from the assumption that $z$ is any point on an arbitrary line $L \subset \mathbb{C} \setminus \{0\}$ that does not pass through the origin, so we could say that $z$ *is* a line that does not pass through the origin. Then, there is a vector $u \in \mathbb{C}$, $|u| = 1$ ($u$ is a unit vector), such that $u \cdot z$ is a vertical line to the right of the origin; algebraically, this vertical line can be expressed as

$$u \cdot z = \frac{1}{c + ce^{i\theta}} = \frac{1}{c + c(\cos\theta + i\sin\theta)}, \quad c > 0, \theta \in (-\pi, \pi). \tag{A1.1}$$

Expression (A1.1) refers to a vertical line, because $\mathrm{Re}(u \cdot z) = \frac{1}{2c}$ is a constant value, while $\mathrm{Im}(u \cdot z) = \frac{-c\sin\theta}{(c+c\cos\theta)^2 + (c\sin\theta)^2}$ is a real-valued continuous function of a real variable $\theta$ that decreases monotonically, with domain $(-\pi, \pi)$ and range $(-\infty, +\infty) = \mathbb{R}$; it can be seen that $\mathrm{Im}(u \cdot z) \to +\infty$ when $\theta \to -\pi^+$, and $\mathrm{Im}(u \cdot z) \to -\infty$ when $\theta \to \pi^-$. $\mathrm{Im}(u \cdot z)$ changes most rapidly near $\theta = \pm\pi$, and it is symmetrical with respect to the origin, as can be seen in Figure A1.

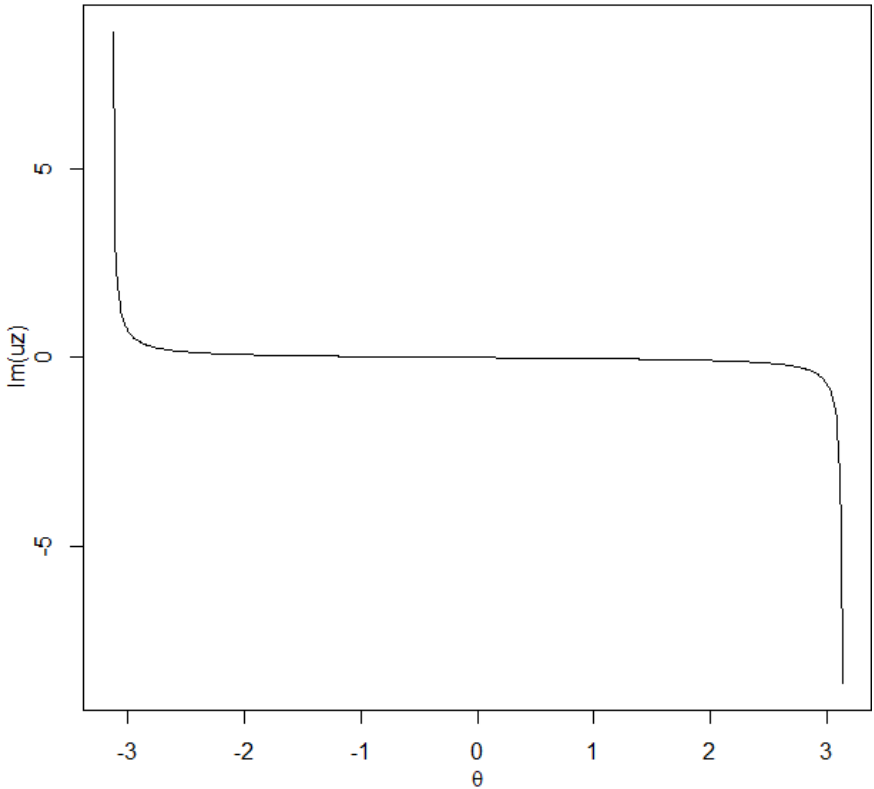

**Figure A1**. Behavior of the imaginary part of $u \cdot z$, within interval $\theta \in (-\pi, \pi)$, for the parametric value $c = 10$. Note that $\mathrm{Im}(u \cdot z) = 0$ when $\theta = 0$.

Now, it is evident that $\frac{1}{u \cdot z} = c + ce^{i\theta}$ is a circumference with center $c$ on the real axis of the plane $\mathbb{C}$, and radius $c$. Then,

$$u \frac{1}{u \cdot z} = \frac{1}{z} = uc + uce^{i\theta} = ce^{i \arg(u)} + ce^{i[\theta + \arg(u)]}. \tag{A1.2}$$

Expression (A1.2) tells us that $\frac{1}{z}$ is a circumference with center $ce^{i \arg(u)}$ and radius $c$. Moreover, circumference $\frac{1}{z}$ passes through the origin (without containing it), since

$$\left| ce^{i \arg(u)} - 0 \right| = \left| ce^{i \arg(u)} \right| = |c| \left| e^{i \arg(u)} \right| = |c| = c; \tag{A1.3}$$

that is to say, the distance between the origin and the center of circumference $\frac{1}{z}$ is equal to its radius. With this, we have proved the lemma stated at the beginning of this demonstration.

To conclude, it is now evident, from expression (A1.2), that if we multiply each point of circumference $1/z$ by a constant value $b \in \mathbb{C}, b \neq 0$, $b/z$ is still a circumference passing through the origin (without containing it), albeit now with center $bce^{i \arg(u)}$ and radius $|b|c$. ∎

**Appendix 2**

**Proposition**. The center $c$ of circumference $C$, which passes through three different points $0, w_1, w_2 \in \mathbb{C}$, is given by

$$c = \frac{-y_1(x_2^2 + y_2^2) + y_2(x_1^2 + y_1^2)}{2x_1 y_2 - 2x_2 y_1} + i \frac{x_1(x_2^2 + y_2^2) - x_2(x_1^2 + y_1^2)}{2x_1 y_2 - 2x_2 y_1} \tag{A2.1}$$

where $x_1 = \text{Re}(w_1)$, $y_1 = \text{Im}(w_1)$, $x_2 = \text{Re}(w_2)$, and $y_2 = \text{Im}(w_2)$. Also, the radius of $C$ is $|c|$.

**Demonstration**: The center $c = h + ik$ and radius $r$ of the circumference with known points $0$, $w_1 = x_1 + iy_1$, $w_2 = x_2 + iy_2$, are obtained by solving the following system of three equations in three unknowns $h, k, r$:

$$(x_1 - h)^2 + (y_1 - k)^2 = r^2$$
$$(x_2 - h)^2 + (y_2 - k)^2 = r^2 \tag{A2.2}$$
$$(0 - h)^2 + (0 - k)^2 = r^2$$

Each of the three equations in system (A2.2) comes from the formula for calculating the squared distance between a known point on a circumference and its center; the third equation in system (A2.2) tells us immediately that $r = |c|$. If we expand the squared binomials in the first two equations of (A2.2) and subtract the third equation of (A2.2) from each of these, we obtain a reduced system of two linear equations in two unknowns $h, k$:

$$2x_1 h + 2y_1 k = x_1^2 + y_1^2 \tag{A2.3}$$
$$2x_2 h + 2y_2 k = x_2^2 + y_2^2$$

The solution of system (A2.3) is

$$h = \frac{-y_1(x_2^2+y_2^2)+y_2(x_1^2+y_1^2)}{2x_1y_2-2x_2y_1} \tag{A2.4}$$

$$k = \frac{x_1(x_2^2+y_2^2)-x_2(x_1^2+y_1^2)}{2x_1y_2-2x_2y_1}$$

The expressions in (A2.4) correspond to the real and imaginary parts on the right side of expression (A2.1). Solution (A2.4) holds as long as $x_1y_2 \neq x_2y_1$ (if $x_1y_2 = x_2y_1$, this means that the points 0, $w_1$, $w_2$ are collinear). ∎

**Appendix 3**

**Proposition**. If line $L_1 \subset \mathbb{C}\backslash\{0\}$, defined in expression (2), contains the roots $r_1, r_2$ of equation $x^2 + c_1x + c_2 = 0$, and $r_1 \neq r_2$, then $\angle 0p_1r_1 = \angle r_1p_1(c_2/p_1)$, with $p_1 = -c_1/2 = (r_1+r_2)/2$.

**Demonstration**: Let's consider the following rectilinear segments:

$A \equiv \overline{0p_1}$,  $\quad\quad B \equiv \overline{p_1r_1}$,  $\quad\quad C \equiv \overline{0r_1}$;

$A' \equiv \overline{r_1p_1}$,  $\quad\quad B' \equiv \overline{p_1(c_2/p_1)}$,  $\quad\quad C' \equiv \overline{r_1(c_2/p_1)}$.

Thus, $\angle 0p_1r_1$ is the angle formed by segments $A$ and $B$, and $\triangle 0p_1r_1$ is the triangle formed by segments $A$, $B$ and $C$; analogously, $\angle r_1p_1(c_2/p_1)$ is formed by segments $A'$ and $B'$, and $\triangle r_1p_1(c_2/p_1)$ is formed by segments $A'$, $B'$ and $C'$. If we denote (for example) the length of the segment $B$ as $|B|$, we can see that $|B| = |\overline{p_1r_1}| = |r_1 - (r_1+r_2)/2| = |(r_1-r_2)/2|$; using this basic approach, it is feasible to verify that

$$\frac{|A|}{|A'|} = \frac{|B|}{|B'|} = \frac{|C|}{|C'|} = \frac{|r_1+r_2|}{|r_1-r_2|}.$$

This last expression tells us that triangles $\triangle 0p_1r_1$ and $\triangle r_1p_1(c_2/p_1)$ are similar, because their corresponding sides are proportional; therefore, the angle formed by segments $A$ and $B$ is equal to the angle formed by corresponding segments $A'$ and $B'$; i.e., $\angle 0p_1r_1 = \angle r_1p_1(c_2/p_1)$. ∎